\documentclass[10pt]{article}
\usepackage{amsmath}
  \usepackage{paralist}
  \usepackage{graphics} 
  \usepackage{epsfig} 
\usepackage{graphicx}  \usepackage{epstopdf}
 \usepackage[colorlinks=true]{hyperref}
 \usepackage{calrsfs}
\usepackage{amssymb}
\usepackage{ifpdf}
\usepackage{graphicx,amssymb,lineno}
\ifpdf
\usepackage[%
  pdftitle={Instructions for use of the document class
    elsart},%
  pdfauthor={Simon Pepping},%
  pdfsubject={The preprint document class elsart},%
  pdfkeywords={instructions for use, elsart, document class},%
  pdfstartview=FitH,%
  bookmarks=true,%
  bookmarksopen=true,%
  breaklinks=true,%
  colorlinks=true,%
  linkcolor=blue,anchorcolor=blue,%
  citecolor=blue,filecolor=blue,%
  menucolor=blue,pagecolor=blue,%
  urlcolor=blue]{hyperref}
\else

\makeatletter
\def\elsartstyle{%
    \def\normalsize{\@setfontsize\normalsize\@xiipt{14.5}}
    \def\small{\@setfontsize\small\@xipt{13.6}}
    \let\footnotesize=\small
    \def\large{\@setfontsize\large\@xivpt{18}}
    \def\Large{\@setfontsize\Large\@xviipt{22}}
    \skip\@mpfootins = 18\p@ \@plus 2\p@
    \normalsize
}
\@ifundefined{square}{}{}
\makeatother

\def\be{\begin{equation}}
\def\ee{\end{equation}}

\begin{document}

\begin{center}{\large{\bf A current value Hamiltonian Approach for Discrete time Optimal Control Problems arising in Economic Growth
Theory}\\

\vspace{1cm}

Rehana Naz }\end{center} { Department of Mathematics and Statistical
Sciences,
  Lahore School of Economics, Lahore, 53200, Pakistan
 Email: drrehana@lahoreschool.edu.pk.\\
   Tel:
 92-3315439545\\
}

\setlength{\parindent}{6mm} \vspace{0.8cm}
\noindent{\bf Abstract} \\
Pontrygin-type maximum principle is extended for the present value
Hamiltonian systems and current value Hamiltonian systems of
nonlinear difference equations for uniform time step $h$. A new
method termed as a discrete time current value Hamiltonian method is
established for the construction of first integrals for current
value Hamiltonian systems of ordinary difference equations arising
in Economic growth theory.
\section{Introduction}
The dynamic optimization problems in economic growth theory are of
fundamental importance for continuous as well as for discrete cases.
The discrete time dynamic optimization problems can be solved by
three methods viz calculus of variation techniques, dynamic
programming and optimal control techniques developed for difference
equations. Halkin \cite{halkin} provided the  maximum principle of
the Pontryagin type \cite{pont1,pont} for systems described by
nonlinear difference equations.  The long-run steady state in a
simple dynamic model of equilibrium with discounted utility function
was studied by Becker \cite{bec}. Day \cite{day} investigated the
irregular growth cycles for the traditional Solow model in discrete
time.  Bohm and Kaas \cite{bohm} and Brianzoni et al \cite{brian}
focused on the discrete time Solow-Swan model with differential
savings. In \cite{brian}, a non-concave production function was
taken into account. Brida and Pereyra  \cite{brida} formulated the
Solow model in discrete time with decreasing population growth rate.
Benhabib et al \cite{ben1} studied the indeterminacy and cycles in
two-sector discrete-time model. A few fundamental discrete time
models in economic were presented by Nishimura and Stachurski
\cite{nis}.

A sperate strand of literature analyzed the discrete time dynamic
optimization problems by dynamic programming.  The necessary
condition for optimality is stated as a Bellman equation \cite{bell}
for the value function and is widely used by economists see e.g
\cite{mer,bell1,bell2}. In the previous literature, the discrete
time economic growth models have been either solved along balanced
growth paths (and the equilibrium) or numerical solutions were
established.

 The first integrals or conservation laws for differential equations are
essential in constructing closed form solutions for both continuous
and discrete cases. The discrete analogue of celebrated Noether's
theorem \cite{Noe} for the discrete variational equations was
provided by Logan \cite{l1} and Dorodnitsyn \cite{l2,l3}. The
discrete Legendre transformation provides the equivalence of the
discrete Euler-Lagrange and discrete Hamiltonian equations
\cite{lal}. Dorodnitsyn and Kozlov \cite{rom1,rom2,vald} established
the relation between symmetries and first integrals for the discrete
Hamiltonian equations by utilizing the discrete Legendre
transformation. The discrete Noether's theorem was formulated in
terms of the discrete Hamiltonian function and symmetry operators.

The optimal control problems involving current value Hamiltonian in
economic growth theory results in a current value Hamiltonian system
which is not of canonical form for both continuous and discrete
cases. Conrad and Clark \cite{con} developed the notion of discrete
time current value Hamiltonian and the maximum principle of the
Pontryagin type for current value Hamiltonian systems described by
nonlinear difference equations.  A current value Hamiltonian
approach was developed by Naz et al \cite{naz,naz2016} to derive the
first integrals and closed-form solutions of current value
Hamiltonian system for continuous case.  The question arises can we
develop a discrete time current value Hamiltonian approach for
difference equations. In this paper,
  the discrete time current value Hamiltonian for
 difference equations arising in economic growth theory is
 developed. For continuous, time a partial Lagrangian approach for dynamical
systems was provided by Naz et al \cite{cl} for models in economics
growth theory. Kara and Mahomed \cite{kar3} developed a partial
Noether approach for PDEs. For approximate ODEs the partial
Lagrangian was proposed by Naeem and Mahomed \cite{nae}.

This paper adds to the literature in multiple ways: firstly the
Pontrygin-type maximum principle is extended for the present value
Hamiltonian systems and current value Hamiltonian systems of
nonlinear difference equations for uniform time step $h$. Secondly,
a current value Hamiltonian approach is developed to establish first
integrals and closed-from solutions of first order ordinary
difference equations arising from sufficient conditions of optimal
control problems in economic growth theory. The current value
Hamiltonian system is not of canonical form so existing techniques
do not hold here.

 The layout of
the paper is as follows. In Section 2, the overview of some basic
notations, definitions and theorems essential for difference
equations is given.  In Section 3, the Pontrygin-type maximum
principle is extended for the present value Hamiltonian systems and
current value Hamiltonian systems of nonlinear difference equations
for uniform time step $h$.  In Section 4, a discrete time
current-value value Hamiltonian approach for construction of first
integrals is proposed.
%

 Finally, conclusions are
summarized in Section 5.

\section{Preliminaries}
I provide here the overview of some basic notations, definitions and
theorem that are essential for this work. The following definitions
and results are adapted from \cite{lal,rom1,rom2, vald}.

\subsection{Hamiltonian formalism for difference equations }
   Consider the space $\bf {\Omega}$ of sequence $(t,\bf
{q}, {\bf p})$.  I introduce one dimensional difference mesh \be
q_{t+1}=q_t+h,\; t=0,1, 2, \cdots, \label{1p}\ee and \be
q_{t-1}=q_t-h,\; t=0,1, 2, \cdots, \label{2p}\ee where a uniform
mesh of step $h$ is considered.

{\it Definition 1:} The Lie operator is defined as following: \be
X=\xi_t\frac{\partial}{\partial t}+\eta^i_{t}
\frac{\partial}{\partial q^i_t}+\zeta^i_{t} \frac{\partial}{\partial
p^i_t}  \label{op}\ee where
\begin{eqnarray} \xi_t=\xi(t,q^i_t,p^i_t),\;
\eta^i_t=\eta^i(t,q^i_t,p^i_t),\;
\zeta^i_t=\zeta^i(t,q^i_t,p^i_t).\nonumber
\end{eqnarray}

{\it Definition 2:} The total shift (left and right) operators and
corresponding discrete operators on a uniform mesh are defined as:
\be S_{+ h}f(t)=f(t+1),\; D_{+ h}=\frac{S_{+ h}-1}{ h},
\label{3p}\ee \be S_{- h}f(t)=f(t-1),\; D_{-h}=\frac{1-S_{- h}}{ h}.
\label{4p}\ee The operators $S_{+h},S_{-h},D_{+h}$ and $D_{-h}$
commute in any combination. Also $D_{+h}=D_{-h}S_{+h}$ and
$D_{-h}=D_{+h}S_{-h}$.

{\it Definition 3:} The discrete Leibniz rule for the operators of
right and left discrete differentiation are as follows: \be D_{+h}(F
G)= D_{+h}(F) G+FD_{+h}( G)+hD_{+h}(F)D_{+h}( G),\ee \be D_{-h}(F
G)= D_{-h}(F) G+FD_{-h}( G)+hD_{-h}(F)D_{-h}( G).\ee

 {\it Definition 4:} The
variational operators are \be \frac{\delta}{\delta
p_t^i}=\frac{\partial}{\partial
p^i_t}+S_{-h}\frac{\partial}{\partial p^i_{t+1}},
i=1,2,\cdots,n,\label{5p}\ee \be \frac{\delta}{\delta
q_t^i}=\frac{\partial}{\partial
q^i_t}+S_{-h}\frac{\partial}{\partial q^i_{t+1}},
i=1,2,\cdots,n,\label{6p}\ee \be \frac{\delta}{\delta
t}=\frac{\partial}{\partial t}+S_{-h}\frac{\partial}{\partial
(t+1)}, \label{7p}\ee where $S_{+ h}$ and $S_{- h}$ are total right
shift and left shift operators as defined in equations {(\ref {3p})}
and {(\ref {4p})}.

 {\it Definition 5:}
 The
action of operators {(\ref {5p})}-{(\ref {6p})} on the
finite-difference functional \be \sum_{\Omega} \bigg(
p^i_{t+1}(q^i_{t+1}-q^i_t)-H(t,t+1,q^i_t,p^i_{t+1})h \bigg),\ee
equated to zero yields following $2n+1$ discrete Hamiltonian
equations:
\begin{eqnarray}
 \frac{q^i_{t+1}-q^i_t}{h}=\frac{\partial H_t}{\partial p^i_{t+1}},
\nonumber\\
 \frac{p^i_{t+1}-p^i_t}{h}=-\frac{\partial H_t}{\partial q^i_t},\;
i=1,\ldots,n, \label{8p}\\
h\frac{\partial H_t}{\partial t}+h\frac{\partial H_{t-1}}{\partial
t}-H_t+H_{t-1}=0,\nonumber
\end{eqnarray}
where \be H_t=H(t,t+1,q^i_t,p^i_{t+1})\:{\text and}\;
H_{t-1}=H(t-1,t,q^i_t,p^i_{t-1})=S_{-h} (H_t).\label{9p}\ee
\subsection{Noether-type theorem for difference Hamiltonian equations}
 The invariance of a difference Hamiltonian on a specified mesh
 provides first integrals of the discrete Hamiltonian equations. To
 consider difference equations three points of mesh are required.

 {\it Definition 6:} The prolongation of the Lie group operator {(\ref {op})}  for neighboring points $(t-1,q_{t-1},p_{t-1})$
 and $(t+1,q_{t+1},p_{t+1})$ is defined as follows:
\begin{eqnarray} X=\xi_t\frac{\partial}{\partial
t}+\eta^i_{t} \frac{\partial}{\partial q^i_t}+\zeta^i_{t}
\frac{\partial}{\partial p^i_t} +\xi_{t+1}\frac{\partial}{\partial
(t+1)}+\eta^i_{t+1} \frac{\partial}{\partial
q^i_{t+1}}+\zeta^i_{t+1} \frac{\partial}{\partial p^i_{t+1}}
\label{10p}\\+\xi_{t-1}\frac{\partial}{\partial (t-1)}+\eta^i_{t-1}
\frac{\partial}{\partial q^i_{t-1}}+\zeta^i_{t-1}
\frac{\partial}{\partial p^i_{t-1}},
\end{eqnarray}
where \begin{eqnarray}
\xi_{t+1}=\xi(t+1,q^i_{t+1},p^i_{t+1}),\;\xi_{t-1}=\xi(t-1,q^i_{t-1},p^i_{t-1}),\nonumber\\
\eta^i_{t+1}=\eta^i(t+1,q^i_{t+1},p^i_{t+1}),\;\eta^i_{t-1}=\eta^i(t-1,q^i_{t-1},p^i_{t-1}),\nonumber\\
\zeta^i_{t+1}=\zeta^i(t+1,q^i_{t+1},p^i_{t+1}),\;\zeta^i_{t-1}=\zeta^i(t-1,q^i_{t-1},p^i_{t-1}).\nonumber
\end{eqnarray}
Dorodnitsyn and Kozlov
 \cite{rom1,rom2} provided following invariance conditions and
 formula for first integrals of the discrete Hamiltonian equations.

{\it Theorem 1:}  A Hamiltonian function is invariant up to a gauge
term $B(t,q_t,p_t)$ with respect to a group generated by operator
{(\ref {op})} if and only if the conditions
\begin{eqnarray}\zeta^i_{t+1} D_{+h}(q_t^i) + p^i_{t+1} D_{+h}(\eta ^i_t)-X(H_t)-H_t D_{+h}(\xi_t)=D_{+h}(B), \label{11p}\\
\nonumber X({\bf \Omega}) |_{{\bf \Omega}=0}=0\end{eqnarray}

hold.

{ \it Theorem 2} (Hamiltonian version of Noether's theorem): The
discrete Hamilton system {(\ref {8p})} which is invariant has the
first integral
 \be I= p^i_t \eta ^i_t -\xi_t( H_{t-1}+h\frac{\partial H_{t-1}}{\partial t})-B \label{12p}\ee
 for some gauge function $B=B(t,q_t,p_t)$ if and only if the Hamiltonian
 action is invariant up to divergence with respect to the operator $X$ given in {(\ref {op})}
on the solutions to equations (\ref{8p}).

{\it Remark 1:} It is important to mention here that
\cite{rom1,rom2, vald} used index free version. I have used slightly
different notions as are used ib \cite{rom1,rom2, vald}. For example
element $q^i_{t+1}$ is denote by $q_i^+$ and similar notions are
used for other variables.
\section{Discrete time Pontryagin type maximum principle and current value Hamiltonian formulation}
In this section, I state the discrete time optimal control problem
of economic growth theory for the infinite horizon for $n$ state,
$n$ costate and $m$ control variables. The discrete time Pontryagin
type maximum principle formulated by \cite{halkin} is presented on a
uniform mesh.  The maximum principle of the Pontryagin type for
current value Hamiltonian systems introduced by Conrad and Clark
\cite{con} is extended for the uniform mesh $h$.
\subsection{Discrete time Pontryagin type maximum principle}
In order to state a discrete time optimal control problem, I
consider the discrete space $ {\bf \Omega}$ for time periods
$t={0,1,2,\cdots, }$. Let $t$ be independent variable, $q_t$ be
vector of $n$ state variables, $ \lambda_t$ be vector of $n$
co-state variables and $u_t$ be the vector of $m$ control variables.
It is worthy to mention here that for discrete case all the
variables are in the form of a sequence for time periods
$t=0,1,\cdots,$.

  In economic analysis, the optimal control
problem is stated as
\begin{eqnarray}
{\rm Maximize}\; {\mathcal F}=\sum_{t=0}^{\infty} \beta ^t F(q_t,u_t) \nonumber\\
{\rm subject\; to}\; \frac{q^i_{t+1}-q^i_t}{h}=f^i(q_t,u_t),\quad
i=1,\ldots,n,\label{1}
\end{eqnarray}
where $\beta<1$ is discount factor.
 Also some appropriate boundary conditions are
imposed. The present value Hamiltonian for discrete time optimal
control problem {(\ref {1})} is defined as \be
H(t,t+1,q_t,\lambda_{t+1},u_t)=\beta ^t
F(q_t,u_t)+\lambda^i_{t+1}f^i(q_t,u_t),\label{2} \ee where
$\lambda^i_{t+1}$ is the co-state variable. The discrete time
Pontryagin type maximum principle is provided by Halkin
\cite{halkin} and requires
\begin{eqnarray}
\frac{\partial H}{\partial u^i_t}=0,\; \frac{\partial^2 H}{\partial u{^i_t}^2}<0\nonumber\\
 \frac{q^i_{t+1}-q^i_t}{h}=\frac{\partial H}{\partial \lambda^i_{t+1}},
\label{3}\\
 \frac{\lambda^i_{t+1}-\lambda^i_t}{h}=-\frac{\partial H}{\partial q^i_t},\;
i=1,\ldots,n.\nonumber
\end{eqnarray}
It is important to mention here that in economic growth theory one
can also have optimal control problems with no discount factor and
for those problems $\beta=1$. Halkin \cite{halkin} provided discrete
time maximum principle for the case $h=1$ and $\beta=1$.
\subsection{Discrete time Current-value Hamiltonian formulation}
 The discrete time optimal control problem {(\ref {1})} is
 non-autonomous but time $t$ enters into picture as part of discount
 factor.  However, it is easy to convert the problem into autonomous
 one by introducing following discrete current value Hamiltonian
 formulation. The discrete time current value Hamiltonian for
 optimal control problem {(\ref {1})} is defined as (see e.g. \cite{con}) \be
 H^c(q_t,p_{t+1},u_t)= F(q_t,u_t)+\beta
p^i_{t+1}f^i(q_t,u_t),\label{4} \ee where \be
H^c=\frac{H}{\beta^t},\;
p^i_{t+1}=\frac{\lambda^i_{t+1}}{\beta^{t+1}}.\label{5}\ee\\
{\bf \it Lemma 1:}  The discrete time Pontryagin type maximum
principle for current value formulation for uniform time step $h$ is
\begin{eqnarray}
\frac{\partial H^c}{\partial u^i_t}=0,\; \frac{\partial^2 H^c}{\partial u{^i_t}^2}<0\nonumber\\
 \frac{q^i_{t+1}-q^i_t}{h}=\frac{\partial H^c}{\partial \beta p^i_{t+1}},
\label{6}\\
 \frac{p^i_{t+1}-p^i_t}{h}=-\frac{\partial H^c}{\partial q^i_t}+\frac{1}{h}(1-\beta)p^i_{t+1},\;
i=1,\ldots,n.  \nonumber
\end{eqnarray}
{\bf \it Proof:\\} The proof of lemma 1 is straight forward. With
the aid of transformations defined in {(\ref {5})} and the discrete
time present value Hamiltonian system {(\ref {4})} one can directly
obtain the discrete time current value Hamiltonian system {(\ref
{6})}.

 The first equation in {(\ref {6})} yields
$u^i_t=g(q_t^i,p^i_{t+1})$ and control variable can be eliminated.

 {\bf
\it Remark 2:} The discount factor $\beta <1$ and the discrete time
current value Hamiltonian system {(\ref {6})} is obtained. Conrad
and Clark \cite{con} provided discrete time maximum principle for
current value Hamiltonian system for the case $h=1$. For $\beta=1$,
the current value Hamiltonian does not exist and optimal control
problems are solved with aid of present value Hamiltonian.

\section{A discrete time current-value value Hamiltonian approach}
The discrete time current-value value Hamiltonian function is
proposed here. The Theorems 1 and 2 do not hold for these as current
value Hamiltonian is not a standard Hamiltonian. An extension of
existing results for this case is needed which I carry out in this
section.

 The discrete time current value Hamiltonian system satisfies

\begin{eqnarray}
 \frac{q^i_{t+1}-q^i_t}{h}=\frac{\partial H}{\partial \beta p^i_{t+1}},
\label{p30}\\
 \frac{p^i_{t+1}-p^i_t}{h}=-\frac{\partial H}{\partial q^i_t}+\Gamma^i,\;
i=1,\ldots,n, \nonumber
\end{eqnarray}
where $\Gamma^i$ are non-zero functions of $t,q^i_t,p^i_{t+1}$. The
extension of the existing results for the standard Hamiltonian
related to discrete time Hamiltonian equations to the discrete time
current-value value Hamiltonian system {(\ref {p30})} is essential
so that one can establish first integrals of system {(\ref {p30})}
in a systematic way.

The current-value Hamiltonian operators determining equations and
formula for first integrals are provided below.

 {\it
Theorem 3:}  A discrete time current value Hamiltonian is invariant
up to a gauge term $B(t,q_t,p_t)$ with respect to a group generated
by operator {(\ref {10p})} if and only if the conditions
\begin{eqnarray}\zeta^i_{t+1} D_{+h}(q_t^i) + p^i_{t+1} D_{+h}(\eta ^i_t)-X(H)-H D_{+h}(\xi_t)\nonumber\\
=D_{+h}(B)+(\eta^i_t-D_{+h}(q^i_t)\xi_t)(-\Gamma^i),
\label{p20}\end{eqnarray} hold, then discrete Hamilton system {(\ref
{8p})} which is invariant has the first integral
 \be I= p^i_t \eta ^i_t -\xi_t H-B. \label{p22}\ee

 {\it Proof:} The discrete derivative of equation {(\ref {p22})}
 with respect to $t$ yields
  \begin{eqnarray} D_{+h}(I)= D_{+h}\bigg(p^i_t \eta ^i_t -\xi_tH-B\bigg)\nonumber\\
 =\eta^i_tD_{+h}(p^i_t)+p_{t+1}^iD_{+h}(\eta^i_t)-HD_{+h}(\xi)-\xi_tD_{+h}(H)-D_{+h}(B)\nonumber\\
 =\eta^i_t (-\frac{\partial H}{\partial q^i_t}+\Gamma^i)+p_{t+1}^iD_{+h}(\eta^i_t)-H D_{+h}(\xi_t)-\xi_t(\Gamma_i\frac{\partial H}{\partial \beta
p^i_{t+1}})-D_{+h}(B)\nonumber\\
=\zeta^i_{t+1} D_{+h}(q_t^i) + p^i_{t+1} D_{+h}(\eta ^i_t)-X(H)-H D_{+h}(\xi_t)\nonumber\\
-D_{+h}(B)-(\eta^i_t-D_{+h}(q^i_t)\xi_t)(-\Gamma^i)\nonumber
  \label{p23},\end{eqnarray}
  where we have used
  \be D_{+h}(H)|_{D_{+h} q^i_t=\frac{\partial H}{\partial \beta p^i_{t+1}},
D_{+h}(p^i_{t})=-\frac{\partial H}{\partial q^i}
+\Gamma_i}=\Gamma_i\frac{\partial H}{\partial \beta
p^i_{t+1}}.\label{a} \ee For $I$ to be first integral  $D_{+h}(I)=0$
and this provides condition as given in {(\ref {p20})}. This
completes the proof.

\section{Conclusions}
Pontrygin-type maximum principle is extended for the present value
Hamiltonian systems and current value Hamiltonian systems of
nonlinear difference equations for uniform time step $h$. A new
method termed as a discrete time current value Hamiltonian method is
established for the construction of first integrals for current
value Hamiltonian systems of ordinary difference equations arising
in Economic growth theory.


\end{document}

